\documentclass[11pt]{article}

\usepackage{makeidx}
\usepackage{amssymb}
\usepackage{amsfonts}
\usepackage{amsmath}
\usepackage{euscript}
\usepackage{theorem}
\usepackage{enumerate}
\usepackage{graphics}

\theorembodyfont{\rm}      

\theoremstyle{change}      

\begingroup

\newtheorem{thm}{Theorem\hskip 5mm}[section]

\newtheorem{cor}[thm]{Corollary\hskip 5mm}

\newtheorem{lem}[thm]{Lemma\hskip 5mm}

\newtheorem{note}[thm]{Note\hskip 5mm}

\endgroup

\begingroup

\endgroup

\begingroup

\endgroup

\newcommand{\qed}{{\unskip\nobreak\hfil\penalty 50\hskip 2em\hbox{}

 \nobreak\hfil$\square$\parfillskip=0pt\finalhyphendemerits=0\par}}

\def\Z{{\bf Z}}

\def\Sp{{\rm Sp}}
\def\GL{{\rm GL}}

\def\l{{\lambda}}

\def\a{{\alpha}}
\def\b{{\beta}}
\def\g{{\gamma}}
\def\w{{\widehat{X}}}
\def\P{{\cal P}}
\def\x{{\cal A}}

\begin{document}

\begin{center}
{\bf Modular reduction of the Steinberg lattice of the general
linear group}
\end{center}

\begin{center}
{\sc Fernando Szechtman}
\end{center}
\begin{center}
{\small\textit{Department of Mathematics \& Statistics, University
of Regina, Saskatchewan, Canada, S4S 0A2}}
\end{center}

\section{Introduction}

We are concerned with a conjecture formulated by Gow [3] regarding
the  reduction modulo $\ell$ of the Steinberg lattice of the
general linear group $G=\GL_n(q)$. Here $n\geq 2$, the underlying
finite field $F_q$ has characteristic $p$, and $\ell$ is a prime
different from $p$.

Let $U$ be subgroup of $G$ consisting of all upper triangular
matrices having 1's along the main diagonal. Write $H$ for the
diagonal subgroup of $G$, and set $B=UH$. Let $\P$ stand for the
set of all subgroups of $G$ that contain $B$, i.e. the standard
parabolic subgroups.

The Steinberg character, say $\chi$, of $G$ is a complex
irreducible character that may be characterized as the only
constituent of the permutation character $1_B^G$ that is not a
constituent of $1_P^G$ for any $P$ in $\P$ different from $B$.
This follows from Steinberg's [4] own determination of $\chi$, and
later work by Curtis [1] expressing $\chi$ as an alternating sum
of permutation characters $1_P^G$.

An explicit realization of $\chi$ was obtained by Steinberg in
[5]. He considers the element $e$ in the integral group algebra
$\Z G$, defined by
$$
e=\left(\underset{\sigma \in
S_n}\sum\mathrm{sign}(\sigma)\sigma\right)\underset{ b\in B} \sum
b,
$$
where the symmetric group $S_n$ is viewed as a subgroup of $G$.
Then the left ideal $I=\Z Ge$ is a $\Z G$-lattice of rank $|U|$
and $\Z$-basis $\{ue\,|\, u\in U\}$, affording $\chi$.
Furthermore, Steinberg shows that  $I/\ell I$ is an irreducible
$F_\ell G$-module if and only if $\ell\nmid [G:B]$. One then is
faced with the interesting problem of trying to find a composition
series for the $G$-module $I/\ell I$ over $F_\ell$, or, rather, a
suitably large extension thereof, when $\ell\mid [G:B]$.

Given that $\chi$ may be realized over $\Z$, restriction to $I$ of
the canonical bilinear form $\Z G\times \Z G\to  \Z$, defined by
$(g,h)\mapsto \delta_{g,h}$, yields a $G$-invariant  symmetric
bilinear form $f:I\times I\to \Z$ with zero radical. Such a form
is unique up to scaling, and Gow normalizes $f$ so that
$f(e,e)=|S_n|$.

We replace the rational integers $\Z$ in the above construction by
a local principal ideal domain $R$ of characteristic 0 and maximal
ideal $\ell R$, containing a primitive $p$-th root of unity. If
$\zeta_p$ is a complex primitive $p$-th root of unity, we may take
$R$ to be the localization of $\Z[\zeta_p]$ at a prime ideal lying
above the unramified prime $\ell$.  The residue field $K=R/\ell R$
has characteristic $\ell$ and a primitive $p$-th root of unity.

The left ideal $I=RGe$ of $RG$ is an $RG$-lattice of rank $|U|$
and $R$-basis $\{ue\,|\, u\in U\}$, affording $\chi$. Note that
$U$ acts on $I$ via the regular representation.

Gow uses the form $f:I\times I\to R$ to produce $RG$-submodules
$I(k)$ of $I$, defined by
$$
I(k)=\{x\in I\,|\, f(x,I)\subseteq \ell^k R\},\qquad k\geq 0.
$$

Consider the $KG$-module $\overline{I}=I/\ell I$ and its
$KG$-submdules $\overline{I(k)}=(I(k)+\ell I)/\ell I$, $k\geq 0$.
This produces the filtration for $\overline{I}$:
\begin{equation}
\label{gowf} \overline{I}=\overline{I(0)}\supseteq
\overline{I(1)}\supseteq \overline{I(2)}\supseteq ...
\end{equation}

As $\overline{I}$ has finite dimension $|U|$ over $K$, only
finitely many of the factors
$M(k)=\overline{I(k)}/\overline{I(k+1)}$ appearing in (\ref{gowf})
are non-zero. By carefully examining the form $f$, Gow is able to
determine the exact non-negative integers $k$ such that $M(k)\neq
(0)$. We will discuss this matter in great detail later.

Based on information from tables for $n\leq 10$, Gow conjectures
that all non-zero $M(k)$ are irreducible $KG$-modules. This would
effectively produce a composition series for $\overline{I}$.

Our contribution to this problem is the following. We start by
showing that the $M(k)$ are completely reducible $KG$-modules.
Next, we prove that the irreducible constituents of a non-zero
$M(k)$ must be self-dual. This is in complete agreement with Gow's
conjecture. It follows that each irreducible constituent of
$\overline{I}$ is self-dual, and we prove that these are pairwise
non-isomorphic. We also show that $\overline{I}$ itself is
self-dual if only if $\overline{I}$ is irreducible.

Our main tool in obtaining these results is the following: if $M$
is a self-dual $KG$-module with no repeated linear characters of
$U$, then $M$ is completely reducible and all its irreducible
constituents are self-dual.

Finally we produce an irreducibility criterion for $M(k)$, and
illustrate its use with a particular case of Gow's conjecture, as
described below.

Let $\kappa_1=\nu_\ell([G:B])$, the $\ell$-valuation of $[G:B]$,
and set $S_1=\overline{I(\kappa_1)}$. It is known that $S_1$ is
irreducible and equal to the socle of $\overline{I}$. Moreover,
all terms
 $\overline{I(s)}$, $s>\kappa_1$,  are known to be $(0)$.
Let $\kappa_2=\nu_\ell([G:P])$, where $P$ is a minimal standard
parabolic subgroup of $G$, and set $S_2=\overline{I(\kappa_2)}$.
Any term of (\ref{gowf}) lying strictly between $S_1$ and $S_2$ is
known to be equal to $S_1$. Moreover, one has $S_2/S_1\neq (0)$ if
and only if $\ell\mid q+1$.

What we show  is that if $\ell$ divides $q+1$ then $S_2/S_1$ is
indeed an irreducible $KG$-module.

Our calculations seem to indicate that $S_2/S_1$ equals the socle
$\overline{I}/S_1$. More generally, for an arbitrary prime $\ell$
dividing $[G:B]$, we wonder if the distinct terms of (\ref{gowf})
are in fact
$$
\overline{I}=S_\omega\supset\cdots \supset S_2\supset S_1\supset
0,
$$
where $S_{i+1}/S_i$ is the socle of $\overline{I}/S_i$.

\section{The underlying root system}

We digress here to develop some notation. Fix a real Euclidean
space with orthonormal basis $e_1,...,e_n$. Then
$$
\Phi=\{e_i-e_j\,|\, 1\leq i\neq j\leq n\}
$$
is a root system in the hyperplane orthogonal to $e_1+\cdots+e_n$.
We use the abbreviated notation
$$
[i,j]=e_i-e_j,\quad 1\leq i\neq j\leq n.
$$
The set
$$
\Pi=\{[i,i+1]\,|\, 1\leq i<n\}
$$
is a fundamental system for $\Phi$, and the associated system of
positive roots is
$$
\Phi^+=\{[i,j]\,|\, 1\leq i<j\leq n\}.
$$
Let $W$ stand for the Weyl group of $\Phi$. We identify $W$ with
the symmetric group $S_n$ via the action of $W$ on
$\{e_1,...,e_n\}$. If $r=[i,j]\in \Phi^+$ then the reflection
$w_r\in W$ is identified with the transposition $(i,j)\in S_n$. As
already mentioned, we view $S_n$, and hence $W$, as a subgroup of
$G$.

If $r=[i,j]\in \Phi$ then $a\mapsto t_r(a)=I+aE^{ij}$ is a group
isomorphism from $F_q^+$ into
$$
X_r=\{I+aE^{ij}\,|\, t\in F_q\}.
$$
Here $I$ stands for the $n\times n$ identity matrix and $E^{ij}$
is the $n\times n$ matrix having a single non-zero entry, namely a
1, in the $ij$-th position. We observe that
\begin{equation}
\label{conj} wt_r(a)w^{-1}=t_{w(r)}(a),\quad w\in W, r\in\Phi,
a\in F_q.
\end{equation}

Suppose now $n\geq 3$ and let $r,s\in\Pi$ be distinct
non-orthogonal roots. Then $r+s\in\Phi^+$ and both $X_r$ and $X_s$
commute elementwise with $X_{r+s}$. For group elements $x,y$ we
denote the commutator $[x,y]=xyx^{-1}y^{-1}$. For $a,b\in F_q$ we
have
\begin{equation}
\label{commu} [t_r(a),t_s(b)]=t_{r+s}((-1)^{\nu(r,s)}ab),
\end{equation}
where $\nu(r,s)=0$ if $r=[i,i+1]$ and $s=[i+1,i+2]$ for some $i$,
while $\nu(r,s)=1$ if $s=[i,i+1]$ and $r=[i+1,i+2]$ for some $i$.

To any subset $J$ of $\Pi$ we associate the standard parabolic
subgroup  $P_J$ of $G$, i.e. the subgroup of $G$ generated by $B$
and all $w_r$, $r\in J$.

\section{A result of Gelfand and Graev}

We will require a modified version of a result due to Gelfand and
Graev, originally proven in the context of complex representations
[2].

This and the following section are the only ones in which the
notation already introduced will be modified. This will allow for
more generality. Here $K$ will stand for an {\em arbitrary} field,
subject solely to the condition of having a primitive $p$-th root
of unity if the characteristic of $K$ is different from $p$. We
will also allow here the case $n=1$.

To stress their dependence on $n$, we let $G(n)=\GL_n(q)$ and
$U(n)=U_n(q)$, the upper unitriangular subgroup of $G(n)$.

We consider the subgroups $H(n)$, $L(n)$ and $A(n)$ of $G(n)$,
defined as follows. For $n\geq 2$ they respectively consists of
all
matrices of the form $$\left(%
\begin{array}{cc}
  1 & u \\
  0 & 1 \\
\end{array}%
\right), \left(%
\begin{array}{cc}
  X & 0 \\
  0 & 1 \\
\end{array}%
\right)\text{ and }\left(%
\begin{array}{cc}
  X & u \\
  0 & 1 \\
\end{array}%
\right),
$$
where $X\in G(n-1)$ and $u$ is a vector in the column space
$F_q^{n-1}$. Note that $A(n)=H(n)\rtimes L(n)$, where $H(n)$ is
canonically isomorphic to $F_q^{n-1}$, and $L(n)$ to $G(n-1)$.
Moreover, we have
\begin{equation}
\label{action}
\left(%
\begin{array}{cc}
  X & 0 \\
  0 & 1 \\
\end{array}%
\right)
\left(%
\begin{array}{cc}
  1 & u \\
  0 & 1 \\
\end{array}%
\right)
\left(%
\begin{array}{cc}
  X & 0 \\
  0 & 1 \\
\end{array}%
\right)^{-1}=\left(%
\begin{array}{cc}
  1 & Xu \\
  0 & 1 \\
\end{array}%
\right).
\end{equation}
We observe that $U(n)$ is a subgroup of $A(n)$. We further define
$H(1)$, $L(1)$ and $A(1)$ to be the trivial subgroups of $G(1)$.
For $n>1$ we will view $A(n-1)$ as canonically embedded in $L(n)$.
Under this embedding $U(n-1)$ becomes a subgroup of $U(n)$, and we
have the decomposition $U(n)=H(n)\rtimes U(n-1)$.

\begin{thm}
\label{gegr}
 Any non-zero module $M$ for $G(n)$ or
$A(n)$ over $K$, whether finite or infinite dimensional, has a one
dimensional subspace that is $U(n)$-invariant.
\end{thm}

\noindent{\it Proof.} The result is clear if $K$ has
characteristic $p$, for in this case the only irreducible
$U(n)$-module is the trivial one. Suppose henceforth that $K$ has
characteristic different from $p$ and that $K$ possesses a
$p$-root of unity different from 1.

Since any $G(n)$-module is automatically an $A(n)$-module, it
suffices to prove the theorem when $M$ is an $A(n)$-module. We
show this by induction on $n$.

The group $U(1)$ being trivial, it acts trivially on $M$, so any
one dimensional subspace of $M$ will do. Suppose next that $n>1$
and the result is true for all $1\leq m<n$. Note that $H(n)$ is a
finite elementary abelian $p$-group. Our assumption on $K$ implies
that
$$
M=\underset{\lambda}\oplus M_\lambda,
$$
where
$$
M_\lambda=\{y\in M\,|\, hy=\lambda(h)y\text{ for all }h\in H(n)\},
$$ and $\lambda$ runs through all group homomorphisms
$H(n)\to K^*$.

Suppose first that $H(n)$ acts trivially on $M$. We consider $M$
as a module for $L(n)$ and, as mentioned above, we view $A(n-1)$
embedded as a subgroup of $L(n)$. By inductive hypothesis there is
a one dimensional subspace of $M$ that is invariant under
$U(n-1)$, and therefore under $U(n)=H(n)U(n-1)$.

Suppose next that $M_\l\neq (0)$ for one or more non-trivial group
homomorphisms $\lambda:H(n)\to K^*$. There is a right action of
$L(n)$ on the set of all non-trivial group homomorphisms
$\mu:H(n)\to K^*$ given by $\mu^x(h)=\mu(xhx^{-1})$ for all $h\in
H(n)$ and $x\in L(n)$. We claim that this action is transitive.
Indeed, our assumption on $K$ yields a non-trivial linear
character $\nu:F_q^+\to K^*$. This gives a non-trivial linear
character $\epsilon:H(n)\to K^*$ defined by
$$
\epsilon\left(%
\begin{array}{cc}
  1 & u \\
  0 & 1 \\
\end{array}%
\right)=\nu(u_{n-1}),\text{ for }u=\left(%
\begin{array}{c}
  u_1 \\
  \vdots \\
  u_{n-1} \\
\end{array}%
\right)\in F_q^{n-1}.
$$
As $\nu$ is non-trivial, formula (\ref{action}) ensures that the
stabilizer of $\epsilon$ in $L(n)$ is exactly $A(n-1)$, again
viewed as subgroup of $L(n)$. Now the index
$[L(n):A(n-1)]=q^{n-1}-1$, which is equal to the number of
non-trivial group homomorphisms $H(n)\to K^*$. This proves the
claim.

Now $M_\lambda\neq (0)$, so $xM_\lambda\neq (0)$ for all $x\in
L(n)$. But clearly $x^{-1}M_\lambda=M_{\l^x}$. By transitivity
$M_{\epsilon}\neq (0)$. Moreover, $M_{\epsilon}$ is invariant
under $A(n-1)$, the stabilizer of $\epsilon$. By inductive
hypothesis there is a one dimensional subspace of $M_{\epsilon}$
that is invariant under $U(n-1)$, and therefore under
$U(n)=H(n)U(n-1)$.\qed

\begin{cor}
\label{equis} Let $M$ be a non-zero $KG$-module. Suppose $M$ has a
one dimensional $U$-invariant subspace $L$ that is contained in
the $KG$-module generated by any other one dimensional
$U$-invariant subspace of $M$. Then the socle of $M$ is
irreducible and is equal to the $KG$-submodule generated by $L$.
Thus $M$ irreducible if and only if it is completely reducible.
\end{cor}

\noindent{\it Proof.} Let $N$ be an irreducible $KG$-submodule of
$M$. From theorem \ref{gegr} we know that $N$ has a one
dimensional $U$-invariant subspace of $M$. By hypothesis $N$
contains $L$. Hence $N$ equals the $KG$-module generated by $L$,
say $S$, which is then irreducible. Thus $M$ has only one
irreducible submodule, namely $S$, so the socle of $M$ is
irreducible and equals $S$.\qed

\section{Complete reducibility of self-dual modules}

As previously mentioned, here we also make modifications to our
general conventions.

\begin{thm}
\label{grande} Let $K$ be a field, $G$ a group, and $U$ a subgroup
of $G$ satisfying: (a) any non-zero $KG$-module has a one
dimensional $U$-invariant subspace; (b) if a $KG$-module admits
$\l:U\to K^*$ as a linear character then it also admits $\l^{-1}$
(for instance, a field $K$ containing a primitive $p$-th root of
unity if the characteristic of $K$ is different from $p$,
$G=\GL_n(q)$ and $U=U_n(q)$, where now $n\geq 1$).

Let $M$ be a self-dual (and hence finite dimensional) $KG$-module.
Suppose that $M$, when viewed as a $KU$-module, has no repeated
irreducible constituents of dimension one. Then $M$ is a
completely reducible (and multiplicity free) $KG$-module.
Moreover, all submodules of $M$ are self-dual as well.
\end{thm}

\noindent{\it Proof.} By assumption there is an isomorphism of
$KG$-modules $\phi:M\to M^*$. To a $KG$-submodule $N$ of $M$ we
associate the $KG$-submodule $N^\perp$ of $M$ defined as follows:
$$
N^\perp=\{x\in M\,|\, \phi(x)(N)=0\}.
$$
Let $P=N\cap N^\perp$ and note that $P\subseteq P^\perp$. We claim
that $P=(0)$. Suppose not. Then by (a) there is a one dimensional
$U$-invariant subspace $L$ of $P$. Then $U$ acts upon $L$ via a
group homomorphism, say $\l:U\to K^*$. By (b) there is a one
dimensional subspace of $P$ upon which $U$ acts via $\l^{-1}$.
This is the same way in which $U$ acts upon $L^*$. Thus $L^*$ is
an irreducible constituent of the $KU$-module $P$, and hence of
$P^\perp$.

Now from $L\subseteq P$ we get the $KU$-epimorphism $P^*\to L^*$.
Likewise, the inclusion $P\subseteq M$ yields the $KG$-epimorphism
$M^*\to P^*$, which can be combined with the $KG$-isomorphism
$M\to M^*$ to produce the $KG$-epimorphism $M\to M^*\to P^*$ with
kernel $P^\perp$. All in all, this yields a $KU$-epimorphism
$M/P^\perp\to L^*$. Hence the multiplicity of $L^*$ as an
irreducible constituent of the $KU$-module $M$ is at least two.
This contradiction proves that $P=(0)$.

Since $N\cap N^\perp=(0)$ and, as noted above, $M/N^\perp\cong
N^*$, we deduce $M=N\oplus N^\perp$. This shows that $M$ is
completely reducible.

The fact that $M$ is multiplicity free follows at once from (a)
and the fact that $M$ has no repeated linear characters of $U$.

To see that $N$ is also self-dual, we consider the linear map
$\phi_N:N\to N^*$ defined by $\phi_N(x)(y)=\phi(x)(y)$ for all
$x,y\in N$. It is a $KG$-homomorphism with kernel $N\cap
N^\perp=(0)$. Since $N$ and $N^*$ have the same dimension,
$\phi_N$ is an isomorphism.

It remains to verify that the given example works. Clearly (a) is
just theorem \ref{gegr}, while (b) can be confirmed as follows.
Let $N$ be $KG$-module with a one dimensional $U$-invariant
subspace $L$ upon which $U$ acts via the group homomorphism
$\l:U\to K^*$. Recall that $H$ stands for the diagonal subgroup of
$G$. Given $h\in H$, the subspace $hL$ of $N$ is also one
dimensional. Furthermore, for $u\in U$ we have
$uhL=hh^{-1}uhL=\l(h^{-1}uh)hL$, so $hL$ is $U$-invariant and is
acted upon by $U$ via a group homomorphism, which we denote by
${}^h \l$. Consider the special element
$h=\mathrm{diag}(-1,1,-1,1,...)$ of $H$. For any $r\in\Pi$ and
$a\in F_q$ we have
$$
h^{-1}t_r(a)h=t_r(-a)=t_r(a)^{-1}.
$$
Since $\l$ is determined by its effect on the fundamental root
subgroups $X_r$, $r\in \Pi$, it follows that ${}^h \l=\l^{-1}$,
the inverse character of $\l$. Thus $U$ acts upon $hL$ via
$\l^{-1}$, as required.\qed

\section{The Steinberg lattice I}

We return to the $RG$-lattice $I$ and its $\ell$-modular reduction
$\overline{I}$. Given $x\in I$, we let $\overline{x}=x+\ell
I\in\overline{I}$. The set $\{u\overline{e}\,|\, u\in U\}$ is a
$K$-basis of $\overline{I}$, so $U$ acts upon $\overline{I}$ via
the regular representation.

To a group homomorphism $\l:U\to R^*$ we associate the set $J(\l)$
of all $r\in\Pi$ such that $\l$ is non-trivial on $X_r$. Let
$P_{J(\l)}$ be the corresponding standard parabolic subgroup of
$G$. We also associate to $\l$ the element $E_\l$ of $RG$ defined
by
$$ E_\l=\underset{u\in U}\sum\l(u)ue.
$$
Clearly $E_\l\neq (0)$ and $xE_\l=\l(-x)E_\l$ for all $x\in U$.
Thus $E_\l$ spans the rank one submodule of $I$ upon which $U$
acts via $\l^{-1}$.

Let $c_\l=\nu_\ell([G:P_{J(\l)}])$. Through skillful calculations,
Gow shows that
\begin{equation}
\label{pre1} E_\l\in I(c_\l),\quad E_\l\notin I(c_\l+1).
\end{equation}
At this point Gow asserts that
\begin{equation}
\label{pre2} \overline{E_\l}\in \overline{I(c_\l)},\quad
\overline{E_\l}\notin \overline{I(c_\l+1)}.
\end{equation}
This crucial fact is true, but it does not follow automatically
from (\ref{pre1}). Indeed, we know from [3] that $f(x,y)=1$ for
some $x,y\in I$, so $x$ belongs to $I(0)$ but not to $I(1)$. Then
$z=\ell x$ belongs to $I(1)$ and not to $I(2)$, but
$\overline{z}=0$, so $\overline{z}$ {\em does} belong to
$\overline{I(2)}$.

A special argument is required to justify (\ref{pre2}). We digress
here to supply the missing details. These are intimately related
to the fact that each non-zero $M(k)$ is a self-dual $KG$-module.

Let  $m=|U|$. As $f$ has zero radical, we see that 0 is not an
elementary divisor of $f$. Thus, these are of the form
$\ell^{a_1},...,\ell^{a_m}$, where the $a_1\leq\cdots\leq a_m$ are
non-negative integers and, by above, $a_1=0$. Let
$\{x_1,...,x_m\}$ and $\{y_1,...,y_m\}$ be $R$-bases of $I$ chosen
so that $f(x_i,y_j)=\ell^{a_i}\delta_{ij}$.

For ease of notation set $c=c_\l$. We identify $\ell^c
R/\ell^{c+1} R$ with $K=R/\ell R$ via the map
\begin{equation}
\label{mapa3} r+\ell R\mapsto \ell^c r+\ell^{c+1} R,\quad r\in R.
\end{equation}

The $R$-bilinear form $f:I\times I\to R$ gives rise to a
well-defined $K$-bilinear form, say $f_c:\overline{I(c)}\times
\overline{I(c)}\to K$, as follows
\begin{equation}
\label{mapa4} f_c(\overline{x},\overline{y})=f(x,y)+\ell^{c+1}
R,\quad x,y\in I(c).
\end{equation}

Clearly $\overline{I(c+1)}$ is contained in the radical of $f_c$.
This naturally produces a $K$-bilinear form
$\overline{f_c}:M(c)\times M(c)\to K$. A basis $B_1$ for $M(c)$ is
formed by all $\overline{x_i}+\overline{I(c+1)}$, if any, such
that $a_i=c$. A basis $B_2$ is obtained by taking the
corresponding $\overline{y_i}+\overline{I(c+1)}$. Taking into
account the identification (\ref{mapa3}) and the definition
(\ref{mapa4}), we see that the matrix of $\overline{f_c}$ relative
the pair of bases $(B_1,B_2)$ is simply the identity matrix of
size $\mathrm{dim}_K M(c)$ (at this point this is possibly zero).
Thus $\overline{f_c}$ is non-degenerate, so the radical of $f_c$
is precisely $\overline{I(c+1)}$. In particular, $M(c)$ is
self-dual.

Why is $M(c)\neq (0)$? Of course, this would follow from Gow's
statement (\ref{pre2}). Why is this statement true? Well, by above
this is equivalent to $\overline{E_\l}$ not being in radical of
$f_c$. Since the group homomorphism $\l^{-1}$, inverse to $\l$, is
associated to the same standard parabolic subgroup as $\l$, it
follows that $E_{\l^{-1}}$ also belongs to $I(c)$. We claim that
$f_c(\overline{E_\l},\overline{E_{\l^{-1}}})\neq 0$, thereby
justifying (\ref{pre2}). Indeed, by virtue of lemma 3.1 and
theorem 3.6 of [3] we have
$$
f_c(\overline{E_\l},\overline{E_{\l^{-1}}})=|U|[G:P_{J(\l)}]+\ell^{c+1}R.
$$
Since $c=\nu_\ell([G:P_{J(\l)}])$ and $\ell\nmid |U|$, our claim
is established.

Taking into account the above discussion and theorem \ref{gegr},
we obtain the following result due to Gow.

\begin{thm}
\label{gow5} Let $k\geq 0$. Then  a factor $M(k)$ of (\ref{gowf})
is not zero if and only if $k=\nu_\ell([G:P])$ for some standard
parabolic subgroup $P$ of $G$. Moreover, in this case:

(a) A linear character $\l:U\to K^*$ enters $M(k)$ if and only if
$\nu_\ell([G:P_{J(\l)}])=k$, in which case it enters only once.

(b) The $KG$-module $M(k)$ is self-dual.
\end{thm}

\section{The Steinberg lattice II}

Here we discuss some of the consequences of the results obtained
in the previous sections. We let $\x$ stand for the set of all
$\nu_\ell([G:P])$ as $P$ ranges through $\P$.

\begin{thm}
\label{exito1} The $KG$-module $\overline{I}$ is multiplicity
free.
\end{thm}

\noindent{\it Proof.} This follows from theorem \ref{gegr} and the
fact that $U$ acts on $\overline{I}$ via the regular
representation, where $\ell\nmid |U|$.\qed

\begin{note} The above result does not hold, in general, for the
$\ell$-modular reduction of the Steinberg lattice of other
classical groups, as no linear character a Sylow $p$-subgroup may
be present in a given composition factor for such a group (cf.
Example 5.4 of [3]). However, we do want to point out that the
multiplicity of the two non-equivalent Weil modules found in [6]
to be constituents of $\overline{I}$ when $\ell=2$ for the
symplectic group $\Sp_{2n}(q)$, $q$ odd, {\em is} indeed one. This
is so because a linear character of the type described above {\em
is} present in each Weil module (cf. section 4 of [6]).
\end{note}

\begin{thm}
\label{exito2} The $KG$-module $\overline{I}$ is self-dual if and
only if it is irreducible.
\end{thm}

\noindent{\it Proof.} If $\overline{I}$ is irreducible then
$\overline{I(1)}=(0)$, so $f_0$ is a non-degenerate $G$-invariant
bilinear form on $\overline{I}$, whence $\overline{I}$ is
self-dual. Conversely if $\overline{I}$ is self-dual, then theorem
\ref{grande} implies that $\overline{I}$ is completely reducible.
But the socle of $\overline{I}$ is known to be irreducible, so
$\overline{I}$ itself is irreducible.\qed

\begin{thm}
\label{exito3}
 Let $k\in \x$. Then $M(k)$ is a self-dual, completely reducible,
non-zero $KG$-module, each of whose irreducible constituents is
also self-dual.
\end{thm}

\noindent{\it Proof.} We know from theorem \ref{gow5} that $M(k)$
is non-zero and self-dual. The remaining assertions follow from
theorem \ref{grande}.\qed

\begin{thm}
\label{exito4} The irreducible constituents of the $KG$-module
$\overline{I}$ are self-dual.
\end{thm}

\noindent{\it Proof.} Taking into account (\ref{gowf}), theorem
\ref{gow5} implies that each irreducible constituent of
$\overline{I}$ must be a constituent of one of the $KG$-modules
$M(k)$, $k\in\x$. By theorem \ref{exito3} each constituent of such
$M(k)$ is self-dual, as required.\qed

\bigskip

We next produce an irreducibility criterion for $M(k)$. Given
$P\in\P$, we see that $H$ acts transitively on the set of group
homomorphisms $\l:U\to K^*$ associated to $P$. We denote by $E_P$
a fixed representative from the $H$-orbit of all $E_\l$, with $\l$
associated to $P$. Note that $E_P$ generates the same
$RG$-submodule of $I$ as any other representative, so whether
$E_P$ belongs to a given $I(k)$ or not depends only on $P$ and not
on the chosen representative.

\begin{thm}
\label{exito5} (Irreducibility Criterion)
 Let $k\in \x$. Then $M(k)$ is an irreducible $KG$-module
if and only if there exists $P\in\P$ such that $\nu_\ell([G:P])=k$
and such that for any other $Q\in\P$ satisfying
$\nu_\ell([G:Q])=k$, the image of $E_P$ in $M(k)$ belongs to the
$KG$-submodule of $M(k)$ generated by the image of $E_Q$.
\end{thm}

\noindent{\it Proof.} Necessity is clear. We know from theorem
\ref{exito3} that $M(k)$ is completely reducible. Thus sufficiency
follows from the above discussion and corollary \ref{equis}.\qed

We wish to apply this criterion to confirm a particular case of
Gow's conjecture. First we need to make sure that the hypotheses
of our criterion are met. This requires three subsidiary results:
one involving the index in $G$ of certain parabolic subgroups, and
two more describing some identities in the group algebra $RG$.

\begin{lem}
\label{mismo}
 Suppose that $\ell$ divides $q+1$.
Let $k=\nu_\ell([G:B])-\nu_\ell(q+1)$. Then $k\in\x$ and the only
standard parabolic subgroups $P$ satisfying $\nu_\ell([G:P])=k$
are those with associated fundamental subset equal to either
$J=\{r\}$, with $r\in \Pi$, or $J=\{r,s\}$, with $r,s\in\Pi$
distinct and non-orthogonal to each other. The last type exists
only if $n\geq 3$.
\end{lem}

\noindent{\it Proof.} Let $J\subseteq \Pi$ and set $Q=P_J$,
$s=\nu_\ell([G:Q])$. If $J=\emptyset$ then $s=k+\nu_\ell(q+1)>k$.
If $J=\{r\}$ for some $r\in\Pi$ then $s=k$. If $J=\{r,s\}$ for
distinct $r,s\in \Pi$ then either $r\not\perp s$, in which case
$[Q:B]=(q+1)(q^2+q+1)$ and $s=k-\nu_\ell(q^2+q+1)=k-0=k$, or else
$r\perp s$, in which case $[Q:B]=(q+1)^2$ and
$s=k-\nu_\ell(q+1)<k$.

If $J\supseteq \{r,s,t\}$ for distinct $r,s,t\in \Pi$ then either
$r,s,t$ are all orthogonal to each other, in which case
$(q+1)^3=[P_{\{r,s,t\}}:B]$ divides $[Q:B]$, or two of them are
non-orthogonal to each other but orthogonal to the third, in which
case $(q+1)^2(q^2+q+1)=[P_{\{r,s,t\}}:B]$ divides $[Q:B]$, or else
one of them is non-orthogonal to the other two, in which case
$(q+1)^2(q^2+q+1)(q^2+1)=[P_{\{r,s,t\}}:B]$ divides $[Q:B]$. In
all  cases $s<k$.\qed

\section{Three identities in $RG$}

Given a group homomorphism $\lambda:F_q^+\to R^*$ and $r\in\Pi$ we
associate to them the group homomorphism $\lambda[r]:U\to R^*$
defined to be trivial on $X_{s}$ for $r\neq s\in\Pi$ and
satisfying $t_{r}(a)\mapsto \lambda(a)$ for all $a\in F_q$.

For $r\in\Pi$, we clearly have
\begin{equation}
\label{hola} w_re=-e.
\end{equation}
Moreover, formula (17) of [5] gives
\begin{equation}
\label{hola2} w_r t_r(a)e=t_r(-a^{-1})e-e,\quad a\neq 0.
\end{equation}
Given a subset $Y$ of $G$ we let $\widehat{Y}$ stand for the
element $\underset{y\in Y}\sum y$ of $RG$. If $r\in \Phi$ we write
$\widehat{X}_r$ rather than $\widehat{X_r}$. From (\ref{hola}) and
(\ref{hola2}) we obtain
\begin{equation}
\label{hola3} w_r \widehat{X}_r e=\widehat{X}_r e -(q+1)e,\quad
r\in \Pi.
\end{equation}
Finally, if $\l:F_q^+\to R^*$ is a group homomorphism we set
$\widehat{X}_{\l(r)}=\underset{a\in F_q}\sum \l(a)t_r(a)$.

\begin{thm}
\label{ex1}
 Assume $n\geq 3$. Let $\lambda:F_q^+\to R^*$ be a non-trivial group
homomorphism. Suppose that $r,s\in\Pi$ are distinct and
non-orthogonal. Then
\begin{equation}
\label{12good} \widehat{X}_{r}w_r\widehat{X}_{s}w_s
E_{\l[r]}=E_{\l[s]}=\widehat{X}_{r}\widehat{X}_{r+s}w_rw_s
E_{\l[r]}.
\end{equation}
Moreover, if $\mu:F_q^+\to R^*$ is a non-trivial group
homomorphism then
\begin{equation}
\label{good}
\widehat{X}_{\mu(r)}w_r\widehat{X}_{s}w_sE_{\l[r]}=(q^2+q+1)
E_{\mu[r]\l[s]}=
\widehat{X}_{\mu(r)}\widehat{X}_{r+s}w_rw_sE_{\l[r]}.
\end{equation}
\end{thm}

\begin{thm}
\label{ex2} Let $n\geq 3$ and let $\lambda,\mu$ be non-trivial
group homomorphisms $F_q^+\to R^*$. Suppose that $r,s\in\Pi$ are
distinct and non-orthogonal. Then
\begin{equation}
\label{hood} \widehat{X}_{r}w_r\widehat{X}_{s}w_s
E_{\l[r]\mu[s]}=E_{\l[s]}=\widehat{X}_{r}\widehat{X}_{r+s}w_rw_s
E_{\l[r]\mu[s]}.
\end{equation}
\end{thm}

It will be better to postpone these two technical proofs until the
end of the paper.

\section{Irreducibility of $S_2/S_1$}

\begin{thm}
\label{casa}
 Suppose $\ell$ divides $q+1$. Let $S_1$ stand for the
socle of $\overline{I}$. Starting at $S_1$ and going up in the
filtration (\ref{gowf}), let $S_2$ be the first term strictly
containing $S_1$. Let $M=S_2/S_1$. Then $M$ is an irreducible
$KG$-module.
\end{thm}

\noindent{\it Proof.} By (\ref{12good}) all $E_{P_J}$, where
$J=\{r\}$ and $r\in\Pi$, generate the same $RG$-submodule of $I$.
By (\ref{hood}) this submodule is contained in the one generated
by any $E_{P_J}$, where $J=\{r,s\}$ and $r,s\in\Pi$ are distinct
and non-orthogonal. By lemma \ref{mismo} the index in $G$ of these
standard parabolic subgroups has the same $\ell$-valuation, which
is shared by no other type of standard parabolic subgroup. Now
apply theorem \ref{exito5} with $P=P_J$, $J=\{r\}$.\qed

\section{Proof of the identities in $RG$}

\noindent{\it Proof of theorem \ref{ex1}.} We will make repeated
and implicit use of (\ref{conj}), (\ref{hola}), (\ref{hola2}) and
(\ref{hola3}) throughout.

The first and third terms (\ref{12good}) are clearly equal, a
comment which also applies to (\ref{good}). We are thus reduced to
proving the second equality in both (\ref{12good}) and
(\ref{good}).

Let $R=U_{w_{r+s}}^+$, i.e. the group of all $u\in U$ such that
$w_{r+s}uw_{r+s}^{-1}\in U$. It equals the product -taken in any
order- of all root subgroups $X_t$, where $t\in \Phi^+$ is
different from $r$, $s$ and $r+s$. We easily see that $R$ is a
normal subgroup of $U$. Thus every element of $X_{r}X_{r+s}$
commutes with $\widehat{R}$, so $\widehat{X}_{r}$ and
$\widehat{X}_{r+s}$ also commute with $\widehat{R}$. Using
(\ref{conj}) we verify that $w_r$ and $w_s$ conjugate
$\widehat{R}$ back into itself. We will implicitly use all these
facts below.

It is easy to see that $U$ acts on the right hand side of
(\ref{12good}) via $\l[s]^{-1}$. Hence this right hand side must
be a scalar multiple of $E_{\l[s]}$. This would be enough for our
purposes, provided the scalar is not 0 modulo $\ell R$. It is not
clear why that should be the case. In fact, a similar remark
applies to (\ref{good}), and in this case the right hand side does
become 0 modulo $\ell R$ when $\ell$ divides $q^2+q+1$. This is to
be expected, since, modulo $\ell I$, $E_{\mu[r]\l[s]}$ lies
strictly above $E_{\l[r]}$ in the filtration (\ref{gowf}) when
$\ell\mid q^2+q+1$.

Returning to the proof, note that, by definition
$$
E_{\l[r]}=\widehat{R}\w_{r+s}\w_{\l(r)}\w_{s}e.
$$
Therefore
\begin{equation}
\label{c1}
w_sE_{\l[r]}=\widehat{R}\w_{r}\w_{\l(r+s)}\w_{s}e-(q+1)\widehat{R}\w_{r}\w_{\l(r+s)}e.
\end{equation}
We first work on the second of the summands appearing on the right
hand side of (\ref{c1}). As $X_{r}$ and $X_{r+s}$ commute
elementwise, we have
\begin{equation}
\label{c2} w_r \w_{r}\w_{\l(r+s)}e =
w_r\w_{\l(r+s)}\w_{r}e=\w_{\l(s)}\w_r e-(q+1)\w_{\l(s)} e.
\end{equation}
Multiplying (\ref{c2}) on the left by
$-(q+1)\widehat{X}_{r+s}\widehat{R}$ we get
\begin{equation}
\label{c3}
\widehat{X}_{r+s}w_r\left(-(q+1)\widehat{R}\w_{r}\w_{\l(r+s)}e\right)=(q+1)^2\widehat{R}\w_{r+s}\w_{\l(s)}e
-(q+1)E_{\l[s]}.
\end{equation}
We now turn our attention to the first summand appearing in
(\ref{c1}). We have
\begin{equation}
\label{c41}
\w_{r}\w_{\l(r+s)}\w_{s}e=\underset{\alpha,\beta,\gamma}\sum
t_{r}(\alpha)
t_{r+s}(\beta)t_{s}(\gamma)\l(\b)e=\underset{\alpha,\beta,\gamma}\sum
 t_{r+s}(\beta)t_{r}(\alpha)t_{s}(\gamma)\l(\b)e,
\end{equation}
since $X_{r+s}$ and $X_r$ commute elementwise. By the commutator
formula (\ref{commu})
\begin{equation}
\label{c42} \underset{\alpha,\beta,\gamma}\sum
 t_{r+s}(\beta)t_{r}(\alpha)t_{s}(\gamma)\l(\b)e=\underset{\alpha,\beta,\gamma}\sum
 t_{r+s}(\beta)t_{s}(\gamma)t_{r}(\alpha)\l(\b+(-1)^{\zeta(r,s)}\a\g)e,
\end{equation}
where $\zeta(r,s)$ is 0 or 1 and depends only on the pair $(r,s)$.
Thus (\ref{c41}) and (\ref{c42}) give
\begin{equation}
\label{c5}w_r \w_{r}\w_{\l(r+s)}\w_{s}e =-\w_{\l(s)}\w_{r+s}e+
\underset{\alpha\neq 0,\beta,\g}\sum
t_{s}(\beta)t_{r+s}(\g)\l(\b+(-1)^{\zeta(r,s)}\a\g)(t_{r}(-\a^{-1})e-e).
\end{equation}
The second summand  in (\ref{c5}) equals
$$
\underset{\alpha\neq 0,\beta,\g}\sum
t_{s}(\beta)t_{r+s}(\g)t_{r}(-\a^{-1})\l(\b+(-1)^{\zeta(r,s)}\a\g)e-
\underset{\alpha\neq 0,\beta,\g}\sum
t_{s}(\beta)t_{r+s}(\g)\l(\b+(-1)^{\zeta(r,s)}\a\g)e.
$$
Using $\underset{\delta}\sum\l(\delta)=0$ we easily verify that
left multiplication of $\widehat{X}_{r+s}$ by each of these two
summands is 0. Therefore (\ref{c5}) gives
\begin{equation}
\label{c6}
\widehat{X}_{r+s}w_r\widehat{R}\w_{r}\w_{\l(r+s)}\w_{s}e=-q\widehat{R}\w_{r+s}\w_{\l(s)}e.
\end{equation}
Going back to (\ref{c1}) and taking into account (\ref{c3}) and
(\ref{c6}) we obtain
\begin{equation}
\label{c7} \widehat{X}_{r+s}w_r
w_sE_{\l[r]}=(q^2+q+1)\widehat{R}\w_{r+s}\w_{\l(s)}e-(q+1)E_{\l[s]}.
\end{equation}
Respectively multiplying (\ref{c6}) on the left by
$\widehat{X}_{r}$ and $\widehat{X}_{\mu(r)}$ yields (\ref{12good})
and (\ref{good}).\qed

\noindent{\it Proof of theorem \ref{ex2}.} As in the previous
result, we only need to prove the second equality. All remarks
made earlier about $R=U_{w_{r+s}}^+$ are still valid. The implicit
use of (\ref{conj}), (\ref{hola}), (\ref{hola2}) and (\ref{hola3})
remains in effect.

By definition
$$
E_{\l[r]\mu[s]}=\widehat{R}\w_{r+s}\w_{\l(r)}\w_{\mu(s)}e.
$$
We let
$$
A=\w_{r+s}\w_{\l(r)}\w_{\mu(s)}e=\underset{\alpha,\beta,\gamma}\sum
t_{r+s}(\alpha) t_{r}(\beta)t_{s}(\gamma)\l(\b)\mu(\g)e.
$$
Then
$$
w_sA=-\w_r\w_{\l(r+s)}e+\underset{\gamma\neq 0,\alpha,\beta}\sum
t_{r}(\alpha) t_{r+s}(\beta)t_{s}(\gamma)\l(\b)\mu(-\g^{-1})e
-\w_r\w_{\l(r+s)}e\underset{\gamma\neq 0}\sum \mu(\gamma).
$$
But
$$
-\underset{\gamma\neq 0}\sum \mu(\gamma)=\mu(0)=1,
$$
so
$$
\begin{aligned}
w_sA &=\underset{\gamma\neq 0,\alpha,\beta}\sum t_{r}(\alpha)
t_{r+s}(\beta)t_{s}(\gamma)\l(\b)\mu(-\g^{-1})e\\
&=\underset{\gamma\neq 0,\alpha,\beta}\sum t_{r+s}(\b)
t_{s}(\gamma)t_{r}(\a)\l(\b+(-1)^{\zeta(r,s)}\a\g)\mu(-\g^{-1})e,
\end{aligned}
$$
where again $\zeta(r,s)$ is 0 or 1 and depends only on the pair
$(r,s)$. Now
$$
w_rw_sA=A_1+A_2+A_3,
$$
where
$$
A_1=-\underset{\gamma\neq 0,\beta}\sum t_{s}(\b)
t_{r+s}(\gamma)\l(\b)\mu(-\g^{-1})e,
$$
$$
A_2=\underset{\gamma\neq 0,\a\neq 0,\beta}\sum t_{s}(\b)
t_{r+s}(\gamma)t_{r}(\a)\l(\b-(-1)^{\zeta(r,s)}\alpha^{-1}\gamma)\mu(-\g^{-1})e,
$$
$$
A_3=-\underset{\gamma\neq 0,\beta}\sum t_{s}(\b)
t_{r+s}(\gamma)\mu(-\g^{-1})a_{\beta,\gamma}e
$$
and
$$
a_{\beta,\gamma}=\underset{\a\neq 0}\sum
\l(\b+(-1)^{\zeta(r,s)}\alpha\gamma).
$$
For $\gamma\neq 0$ and any $\beta\in F_q$ we have
$$
-a_{\beta,\gamma}=\l(\b),
$$
so $A_1+A_3=0$. Hence
$$
w_rw_sA=A_2.
$$
Thus
$$
\w_{r+s}w_rw_sA=\underset{\delta,\gamma\neq 0,\a\neq 0,\beta}\sum
t_{s}(\b)
t_{r+s}(\gamma+\delta)t_{r}(\a)\l(\b-(-1)^{\zeta(r,s)}\alpha^{-1}\gamma)\mu(-\g^{-1})e.
$$
Making the change of variable $\epsilon=\delta+\gamma$ we obtain
$$
\w_{r+s}w_rw_sA =\underset{\epsilon,\a\neq 0,\beta}\sum t_{s}(\b)
t_{r+s}(\epsilon)t_{r}(\a)b_{\a,\b}e,
$$
where
$$
b_{\a,\b}= \underset{\gamma\neq 0}\sum
\l(\b-(-1)^{\zeta(r,s)}\alpha^{-1}\gamma)\mu(-\g^{-1}).
$$
We may write this in the form
$$
\begin{aligned}
\w_{r+s}w_rw_sA &=\w_{r+s}\underset{\a\neq 0,\beta}\sum t_{s}(\b)
t_{r}(\a)b_{\a,\b} e\\
&=\w_{r+s}\underset{\a\neq 0,\beta}\sum t_{r}(\a)t_{s}(\b)
b_{\a,\b} e.
\end{aligned}
$$
The last equality holds because $\w_{r+s}$ absorbs all commutators
arising from elements of $X_{r}$ and $X_{s}$. Multiplying by
$\w_r$ and making a suitable change of variable we get
$$
\w_r\w_{r+s}w_rw_sA=\w_{r+s}\underset{\nu,\beta}\sum t_{r}(\nu)
t_{s}(\b) c_\beta e,
$$
where
$$
c_\beta=\underset{\gamma\neq 0}\sum \mu(-\g^{-1})d_{\beta,\gamma}
$$
and
$$
d_{\beta,\gamma}=\underset{\a\neq 0}\sum
\l(\b-(-1)^{\zeta(r,s)}\alpha^{-1}\gamma).
$$
For $\gamma\neq 0$ and any $\beta$ we have
$$
d_{\beta,\gamma}=-\l(\b).
$$
Therefore
$$
c_\beta=-\l(\b)\underset{\gamma\neq 0}\sum \mu(-\g^{-1})=\l(\b).
$$
Hence,
$$
\w_r\w_{r+s}w_rw_sA=\w_{r+s}\underset{\nu,\beta}\sum t_{r}(\nu)
t_{s}(\b)\l(\b) e=\w_{r+s}\w_{r}\w_{\l(s)}e.
$$
Recalling the meaning of $A$ and multiplying through by
$\widehat{R}$ we obtain (\ref{hood}).\qed

\medskip

\noindent{\bf{\Large Acknowledgments}}

\medskip

I am very grateful to R. Gow for his help in connection with
theorem \ref{gegr}.

\medskip

\noindent{\bf{\Large References}}

\medskip

 \small


\noindent [1] C. W. Curtis, \emph{The Steinberg character of a
finite group with a $(B,N)$-pair}, J. Algebra {\bf 4} (1966),
433-441.

\noindent [2] I.M. Gelfand and M.I. Graev, \emph{Construction of
irreducible representations of simple algebraic groups over a
finite field}, Dokl. Akad. Nauk SSSR, {\bf 147} (1962) 529-532.

\noindent [3] R. Gow, \emph{The Steinberg lattice of a finite
Chevalley group and its modular reduction}, J. London Math. Soc.
(2), {\bf 67} (2003), 593-608.


\noindent [4] R. Steinberg, \emph{A geometric approach to the
representations of the full linear group over a Galois field},
Trans. Amer. Math. Soc. {\bf 71} (1951), 274-282.

\noindent [5] R. Steinberg, \emph{Prime power representations of
finite linear groups II}, Canad. J. Math. {\bf 9} (1957), 347-351.

\noindent [6] F. Szechtman, \emph{The 2-modular reduction of the
Steinberg representation of a finite Chevalley group of type $C\sb
n$}, J. Group Theory {\bf 8} (2005), 11--41.

\end{document}